\documentstyle[11pt,epsf,psfrag]{article}
\input amssym.def
\input amssym.tex
\font\bbb=msbm10 scaled 1100

\newtheorem{thm}{Theorem}[section]
\newtheorem{lemma}[thm]{Lemma}
\newtheorem{cor}[thm]{Corollary}
\newtheorem{prop}[thm]{Proposition}
\newtheorem{conj}[thm]{Conjecture}

\newcommand{\eps}{\epsilon}
\newcommand{\ga}{\gamma}
\newcommand{\del}{\partial}

\newcommand{\ra}{\rightarrow}
\newcommand{\hra}{\hookrightarrow}
\newcommand{\real}{\mbox{\bbb R}} 	% {\mbox{\rm I} \hspace{-.03in} {\bf R}}
\newcommand{\zed}{\mbox{\bbb Z}}	%{\mbox{\bf Z}\hspace{-.06in}{\bf Z}}
\newcommand{\abs}[1]{\vert #1\vert}

\newcommand{\be}{\begin{equation}}
\newcommand{\ee}{\end{equation}}
\newcommand{\bea}{\begin{eqnarray}}
\newcommand{\eea}{\end{eqnarray}}
\newcommand{\bmini}{\footnotesize\begin{center}\begin{minipage}{5.5in}}
\newcommand{\emini}{\end{minipage}\end{center}\normalsize}
\newcommand{\pf}{{\em Proof: }}

\newcommand{\qed}{\hfill$\Box$}

\newcommand{\bfx}{\mbox{{\bf x}}}

\newcommand{\eg}{{\em e.g.}}
\newcommand{\ie}{{\em i.e.}}
\newcommand{\cf}{{\em cf. }}
\newcommand{\goal}{g}
\newcommand{\Ga}{\Gamma}
\newcommand{\Si}{\Sigma}
\newcommand{\Up}{\Upsilon}
\newcommand{\Cs}{{\cal C}}
\newcommand{\Os}{{\cal O}}
\begin{document}

\begin{center}
\large 
{\bf  CONFIGURATION SPACES AND BRAID GROUPS ON GRAPHS IN
	ROBOTICS}
\normalsize
\vspace{0.1in}

%
% AUTHORS
%

ROBERT GHRIST \\ 
School of Mathematics \\
Georgia Institute of Technology \\
Atlanta, GA 30332  

\vspace{0.15in}

{\em Dedicated to Joan Birman on her $70^{th}$ birthday.}

\end{center}
% 
% ABSTRACT
%
\begin{abstract}
Configuration spaces of distinct labeled points on the plane 
are of practical relevance in designing safe 
control schemes for Automated Guided Vehicles (robots) in 
industrial settings. In this announcement, we consider the problem of 
the construction and classification of configuration spaces for 
graphs. Topological data associated to these spaces (\eg, 
dimension, braid groups) provide an effective measure of the 
complexity of the control problem. The spaces are themselves 
topologically interesting objects. We show that they are 
$K(\pi_1,1)$ spaces whose homological dimension is bounded 
by the number of essential vertices. Hence, the braid groups 
are torsion-free.
\vspace{0.1in}

%
% AMS NUMBERS
%

%\noindent
%{\sc keywords: braids, configuration space, robotics.}

\end{abstract}

% &&&&&&&&&&&&&&&&&&&&&&&&&&&&&&&&&&
\section{Configuration Spaces in Manufacturing}
% &&&&&&&&&&&&&&&&&&&&&&&&&&&&&&&&&&

% &-&-&-&-&-&-&-&-&-&-&-&-&-&-&-&-&-&
\subsection{Background}
% &-&-&-&-&-&-&-&-&-&-&-&-&-&-&-&-&-

In several manufacturing and industrial settings, the following 
scenario arises: there is a collection of independent mobile 
Automated Guided Vehicles ({\em aka} AGVs) which traverse a factory 
floor replete with obstacles en route to a goal position (say, a 
loading dock or an assembly workstation), from whence the process
iterates. For these applications, it is of the utmost importance 
to design a control scheme which insures that 
(1) the AGVs not collide with the obstacles; (2) the AGVs
not collide with each other; (3) the AGVs complete the assigned
task with a certain efficiency with respect to various work 
functionals. 

In practice, control schemes are often effected by employing 
high factor-of-safety algorithms which guarantee safety but reduce
efficiency. For example, one partitions the factory floor into 
``zones'' and then all algorithms are written such that only one 
AGV is allowed in a zone at any given time 
\cite{agvhandbook,bozer.iie91}.
Clearly, such practices are not optimally efficient. 

Within the past decade, the idea of using abstract configuration 
spaces to model the workspace and then fabricating a control 
scheme on this topological space, has filtered through the 
robotics community \cite{latombe.book,kod&rimon.aam}. 
Surprisingly, topologists have been generally unaware of and uninvolved 
in this important development.

The technique is very straightforward: assume that the individual 
AGVs are represented as points on the workspace floor $X$. The 
set $\Os\subset X$ represents those obstacles which are
to be avoided. The configuration space $\Cs$ of $N$ noncolliding 
labeled AGVs in the workspace is thus:
\be
\Cs := \left[(X-\Os)\times(X-\Os)\cdots\times(X-\Os)\right]-\Delta,
\ee
where $\Delta$ denotes the pairwise diagonal 
$\Delta := \{(x_1,x_2,\ldots,x_N : x_i=x_j {\mbox{ for some }}
i\neq j\}$. A point in $\Cs$ thus represents a ``safe'' configuration
of AGVs. 

Let $\goal\in\Cs$ denote a desired goal configuration. In order to 
enact a control scheme which realizes this goal safely, it is 
sufficient to build a vector field $X_\goal$ on $\Cs$ which 
(1) has $\goal$ as a sink with a large basin of attraction; and (2) is 
transversally inward on the boundaries left from removing 
the diagonal and the obstacles. By evolving initial states with 
respect to this control field, even initial configurations which are
near a collision are immediately repulsed onto a `safe' pathway. 

This topological/dynamical formulation has several advantages:
\begin{enumerate}
\item
	There are no {\em ad hoc} restrictions on how many 
	AGVs can inhabit a subdomain of the workspace;
\item
	There is an analytical measure of {\em safety} --- the 
	distance to the boundary of $\Cs$ in the natural (product) 
	metric inherited from the workspace --- which allows
	for rigorous treatment of this issue; and
\item
	The control field can be designed to be {\em stable} in the
	sense that, when the basin of attraction of $\goal$ is very 
	large (it often can be a submanifold of full measure), 
	arbitrarily large perturbations in the state of the system 
	(\eg, momentary failure of a steering component, 
	slippage, the effect of debris on the floor, etc.) do not
	affect the attainment of the goal state.
\end{enumerate}

The robotics community has effectively employed these ideas 
into real control schemes which work in certain industrial settings: 
see \cite{latombe.book} for an introduction to and description of 
various techniques. 

It is by no means true, however, that this problem is completely
solved. In the case where the configuration space is a manifold, 
the existence of a vector field which realizes the desired goal and
is repulsive along the unsafe boundaries follows from a simple
application of Morse theory (\cf \cite{kod&rimon.aam}). The 
specification of such a vector field is of course a more difficult
problem which requires detailed knowledge of the configuration 
space (at least on the level of charts), though general formulae are
available for some specific instances \cite{kod&rimon.aam}.

% &-&-&-&-&-&-&-&-&-&-&-&-&-&-&-&-&-&
\subsection{Graphs and Guidepath Networks}
% &-&-&-&-&-&-&-&-&-&-&-&-&-&-&-&-&-&

In this paper, we initiate the use of configuration space methods 
in an industrial setting which has heretofore resisted analysis.

The setting we describe above, in which the automated guided 
vehicles have a full two degree-of-freedom autonomy of planar 
motion, assumes a certain level of sophistication in the machinery. 
A more elementary (and hence, inexpensive and easier to install 
and maintain) system involves AGVs which are constrained to 
a network of guidepath wires, embedded either in the floor
as tracks or suspended from the ceiling as wires 
\cite{agvhandbook}. Such systems are currently in use 
in many industrial settings.

The problem of maneuvering AGV's on a graph is completely 
different from that of a full two-degree-of-freedom planar 
system. On the plane, two AGV's on a collision course may avoid
disaster by ``swerving'' around each other at the last minute:
not so on a 1-d graph. Here, the problem is much more global ---
the control scheme requires information about the global 
structure of the graph as a local perturbation is usually 
insufficient. Hence, it would appear that a topological 
approach could be of the greatest possible benefit. 

\begin{figure}[htb]
\begin{center}
\begin{psfrags}
\psfrag{0}[][]{$v_0$}
\psfrag{1}[][]{$v_1$}
\psfrag{2}[][]{$v_2$}
\psfrag{3}[][]{$v_3$}
\psfrag{a}[][]{$e_1$}
\psfrag{b}[][]{$e_2$}
\psfrag{c}[][]{$e_3$}
	\epsfxsize=1.8in\leavevmode\epsfbox{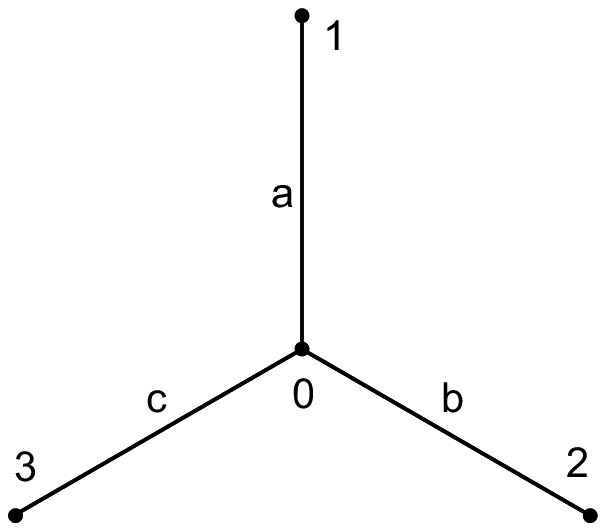}
\hspace{0.5in}
	\epsfxsize=2.9in\leavevmode\epsfbox{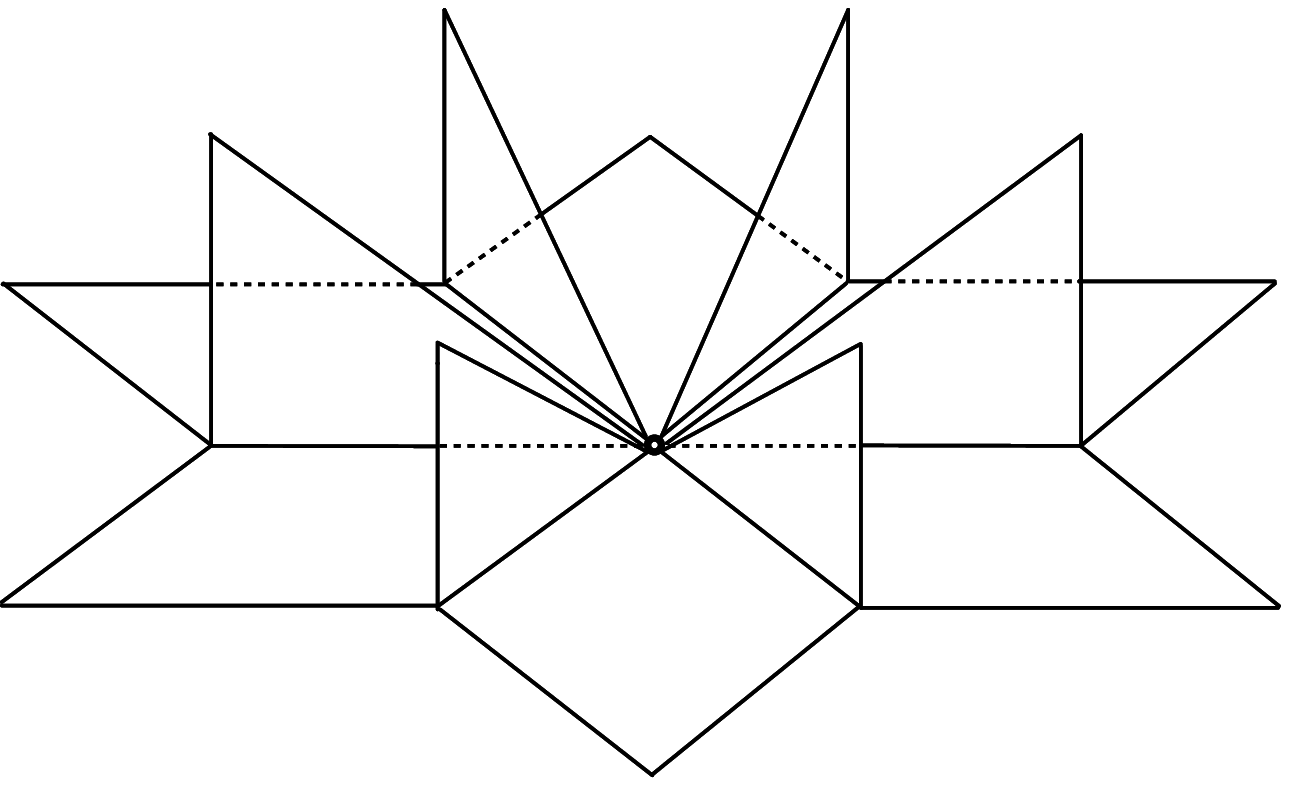}
\end{psfrags}
\end{center}
\caption{The Y-graph $\Ga_Y$ [left] has a configuration space 
	$\Cs^2(\Ga_Y)$ [right] which embeds in $\real^3$.} 
\label{fig_Y}
\end{figure}

We commence our investigation with the simplest nontrivial 
example: that of a pair of points on a `Y-graph' --- the tree
$\Ga_Y$ having three edges $\{e_i\}_1^3$ meeting at a 
single 3-valent vertex $v_0$ as in Figure~\ref{fig_Y}[left]. 
Let $\Cs^2(\Ga_Y)$ denote the configuration space of two points 
on $\Ga_Y$. 

% =========================================
% =========================================

\begin{lemma}
The configuration space $\Cs^2(\Ga_Y)$ is homeomorphic to the 
embedded 2-complex illustrated in Figure~\ref{fig_Y}[right]. 
\end{lemma}
\pf
Let $x$ and $y$ denote the distinct points on $\Ga_y$ and let 
$D\subset\Cs^2(\Ga_Y)$ denote the region where $x$ and $y$ 
lie on distinct edges of $\Ga_Y$. Then $\overline{D}$, the closure
of $D$ within $\Cs^2(\Ga_Y)$ is easily seen to be a 2-manifold
with boundary as follows: if $x$ and $y$ lie on distinct edges, 
the point $(x,y)$ has a product neighborhood; if (say) $x$ lies
on the central vertex, then $y$ lies within one of the three 
edges, and a neighborhood within $\overline{D}$ allows $x$ 
to move onto the other two edges, these two edges together
being homeomorphic to an interval. As there are two points and
three intervals, $\overline{D}$ has a natural decomposition into 
six cells, each homeomorphic to $([0,1]\times[0,1])\setminus\{(0,0)\}$
and joined along the edges $(0,1]\times\{0\}$ and $\{0\}\times(0,1]$
in a cyclic fashion. An Euler characteristic computation 
combined with keeping track of the boundary reveals that 
$\overline{D}$ is homeomorphic to a punctured disc. The complement
of $\overline{D}$ in $\Cs^2(\Ga_Y)$ consists of six ``fins,''
each homeomorphic to $\{(x,y)\in\real^2 : x\geq 0, y\geq 0, x+y<1\}$.
These fins are attached to $\overline{D}$ along the edges 
$\{0\}\times(0,1]$ and $(0,1]\times\{0\}$ in each of the cells.
\qed

Note that this space has a product decomposition as 
the Cartesian product of $(0,1]$ with the graph given by attaching 
six radial edges to a circle. The $(0,1]$ factor is a height function 
on $\Cs^2(\Ga)$ which measures the distance between $x$ and
$y$.

It is worthwhile to contemplate this example for hints as to what 
the general case of $N$ points on a tree might reveal. Note that 
the braid group of this graph is especially simple ($\zed$) and
that the space itself deformation retracts onto a graph. 
For an analysis of the dynamics and control of vector fields on 
this space, see \cite{GK97,GK98}. 

After deriving our main results in the next sections, we conclude
with several additional examples. The interested reader may wish 
to reference these as necessary for intuition.

% &&&&&&&&&&&&&&&&&&&&&&&&&&&&&&&&&&
\section{Configuration spaces of trees}
% &&&&&&&&&&&&&&&&&&&&&&&&&&&&&&&&&&

For the entirety of this section, $\Up$ will denote a 
tree having distinguished vertex $p$. Let $V$ denote the number of 
essential vertices in $\Up$ --- that is, the number of vertices of 
valency greater than two. The configuration space of 
$N$ distinct points on $\Up$ is denoted $\Cs^N(\Up)$. This space
can be analyzed by considering the codimension-one subset
\be
	\Si := \{ \bfx\in\Cs^N(\Up) : x_n=p {\mbox{ for some }} n \} ,
\ee
which splits into $N$ disjoint (and disconnected) components
\be
	\Si_n := \{ \bfx\in\Cs^N(\Up) : x_n=p \} .
\ee
Assume that at $p$ there are $K>2$ incident edges in $\Up$. Then
\be\label{eq_Sigma}
\Sigma \cong \coprod_{n=1}^N\left(\prod_{\abs{j}=N-1}
	\Cs^{j_1}(\Up_1)\times\cdots\times\Cs^{j_K}(\Up_K)\right) ,
\ee
where the $\Up_i$ are the connected components of $\Up-\{p\}$.
The complement of $\Sigma$ has the form 
\be\label{eq_Complement}
	\Cs^N(\Up)-\Si \cong \coprod\left(\prod_{\abs{j}=N}
	\Cs^{j_1}(\Up_1)\times\cdots\times\Cs^{j_K}(\Up_K)\right) .
\ee

Using the subsets $\Si_n$ to decompose $\Cs^N(\Up)$ provides the 
basis for the induction arguments which fill the remainder of 
this announcement.

Consider the case where $\Up$ is a tree having a vertex $p$ 
on its boundary. Let $e$ denote the unique edge connecting $p$ 
to an essential vertex $q$ of $\Up$ and denote by $\Up'$ the 
subtree given by $\Up' := \overline{\Up-\overline{e}}$. 
That is, $\Up'$ contains the point $q$, but neither $e$ nor $p$.
Next define $\Si'_n$ to be the  
set of configurations on $\Up'$ which have the point $x_n$ at $q$;
hence,  $\Si'_n:= \{\bfx\in\Cs^N(\Up') : x_n=q \}$.  Define 
also the ``end'', $\Sigma_n := \{\bfx\in\Cs^N(\Up) : x_n=p \}$.
Note that $\Sigma_n$ is homeomorphic to $\Cs^{N-1}(\Up)$.

\begin{lemma}\label{lem_Decomp}
The space $\Cs^N(\Up)$ is homeomorphic to  
%%%%%%%%%%%%%%%%%%%%%%%%%%%%%%%%%%%%%%%%%%
\be \label{eq_Decomp}
C^N(\Up) \cong C^N(\Up')\bigcup_{\Sigma'_n}^{n=1..N}
\left\{(\Sigma_n\times(0,1])\cup(\Sigma'_n\times\{0\})\right\} ,
\ee
%%%%%%%%%%%%%%%%%%%%%%%%%%%%%%%%%%%%%%%%%%
where the spaces $(\Sigma_n\times(0,1])\cup(\Si'_n\times\{0\})$ 
are glued to $\Cs^N(\Up')$ along $\Si'_n\times\{0\}$.
\end{lemma}
\pf
The key observation is the
following: one can decompose $\Cs^N(\Up)$ into those configurations
in which all points lie on $\Up'$ and those configurations 
in which the point $x_n$ lies on $e$ and is the farthest such point
from $q$: \ie, no point lies on the interval of $e$ from $x_n$ to $p$. 
The set of configurations in which $x_n$ 
is at a fixed point on $e$ (and the farthest such point from $q$) 
is homeomorphic to a copy of $\Sigma_n$, parameterized
by the distance along $e$. This set $\Sigma_n\times(0,1]$ is glued 
to $\Cs^N(\Up')$ along the subset of $\Si'_n\times\{0\}$, 
since no other point may occupy the vertex $q$.
\qed

Using this decomposition, we proceed with the following 
fundamental step.

\begin{lemma}
\label{lem_Injectivity}
For any tree $\Up$ and any vertex $p$, 
the inclusion of $\Sigma$ into $\Cs^N(\Up)$ is 
$\pi_1$-injective.\footnote{Note, everything is basepoint 
dependent since $\Sigma$ is not path connected. This theorem holds
for arbitrary choice of basepoint.}
\end{lemma}
\pf
We induct on the number of essential vertices $V$, letting 
$\Si$ denote the codimension-one subcomplex of configurations
for which there is a point at the $V^{th}$ vertex $p$. For $V=1$, 
the inclusion map is $\pi_1$-injective since every connected
component of $\Sigma$ is contractible. 

Let $\ga$ denote a representative of $\pi_1(\Si)$ which bounds 
a disc $D\subset\Cs^N(\Up)$. By taking $D$ in general position
with respect to $\Si$ and ordering the connected 
components of $\Si\cap D$ with respect
to inclusion, we may assume without a loss of generality that the 
interior of $D$ lies in a connected component of the 
complement of $\Si$. From Equation~(\ref{eq_Complement}) 
it follows that $D$ lies in the product of configuration 
spaces of graphs with strictly fewer essential vertices. 
However, since $D\cap\Si=\ga$, the point $x_n$ is located at 
$p$ along $\ga$, and no other points ever occupy $\Si$ on $D$; 
hence, $x_n$ moves within some subtree of $\Up-\{p\}$.
Thus, it suffices to consider the specialized case where 
$\Up$ is a tree and $p$ is a vertex on the boundary of the tree.
Let $\Si_n$ denote the set of configurations of $N$ distinct 
points on $\Up$ for which the point $x_n$ is at $p$. Assume 
that $\ga$ is a representative of $\pi_1(\Si_n)$ which bounds 
a disc $D\subset\Cs^N(\Up)$. 

Let $e$ denote the unique edge which connects an essential 
vertex $q\in\Up$ to the boundary point $p$. 
Denote by $\Up'$ the subgraph $\overline{\Up-\overline{e}}$. 
Then define the subset $\Si_n' := \{\bfx\in\Cs^N(\Up') : x_n=q \}$.
From Lemma~\ref{lem_Decomp} we have
\be
\Cs^N(\Up) \cong \Cs^N(\Up')\bigcup_{\Si_n'}^{n=1..N}
	\left\{(\Si_n\times(0,1])\cup(\Si'_n\times\{0\})\right\} ,
\ee
where the gluing is along the subset $\Sigma_n'\subset\Si_n\times\{0\}$
at which the $n^{th}$ point is distance ``zero'' from $q$.

We now have $\ga$ a loop in $\Si_n\times\{1\}$  which bounds
a disc $D$ within $\Cs^N(\Up)$. Since $\ga$ is assumed nontrivial 
in $\pi_1$, $D$ must intersect the gluing set $\Sigma_n'$ in 
a nontrivial loop $\ga'\subset\Si_n'$ which bounds a disc $D'$ 
in $\Cs^N(\Up)$. By the induction hypothesis, the inclusion 
$\iota':\Si_n'\ra\Cs^N(\Up)$ is $\pi_1$-injective. Thus, $\ga$ 
is not contractible in $\Cs^N(\Up)$.
\qed

\begin{thm}
\label{thm_Contract}
Given the configuration space $\Cs^N(\Up)$ of a tree $\Up$ 
and a connected subset $K\subset\Cs^N(\Up)$, if the 
homomorphism $\iota_*:\pi_1(K)\ra\pi_1(\Cs^N(\Up))$ 
induced by inclusion is trivial, then $K$ is nullhomotopic
in $\Cs^N(\Up)$. 
\end{thm}
\pf
Induct on the number of essential vertices $V$ of $\Up$.
The theorem is certainly true for $V=0$. As before, let $\Si$ 
denote the configurations which have a point at the $V^{th}$
essential vertex. The complement of $\Si$ in $\Cs^N(\Up)$ 
is composed of products of configuration spaces of graphs 
with fewer numbers of essential vertices, \cf 
Equation~(\ref{eq_Complement});  hence, if $K$ lies within 
the complement of $\Si$, then $K$ contracts by induction. 

In the case where $K\cap\Sigma\neq\emptyset$, we know 
from Lemma~\ref{lem_Injectivity} that $\pi_1(K\cap\Si)\mapsto 0$
under inclusion into $\pi_1(\Si)$. As $\Si$ is composed of products
of configuration spaces of graphs with smaller numbers of essential
vertices (see Equation~(\ref{eq_Sigma})), we have by induction 
that each connected component of $K\cap\Si$ contracts to a point
within $\Si$. One may then ``pinch'' $K$ off of $\Si$ 
into a collection of connected sets $K_i$ since $\Si$ is a codimension
one subcomplex. By induction, we may homotope each $K_i$ 
to a point within the complement of $\Si$. By tracing out the 
path of the pinch-points under the homotopies, we have a homotopy
of all of $K$ to a graph, which must be 
nullhomotopic in $\Cs^N(\Up)$.
\qed

\begin{cor}\label{cor_EM}
The configuration space $\Cs^N(\Up)$ is an Eilenberg-MacLaine 
space of type $K(\pi_1,1)$: \ie, $\pi_k(\Cs^N(\Up))=0$ for all 
$k>1$.
\end{cor}

Corollary~\ref{cor_EM} is significant in that $\pi_1$ 
determines the homotopy type of the configuration 
space. We thus consider the (pure) braid groups of trees.
For $\Up$ a planar graph, the inclusion $\iota:\Up\hra\real^2$ 
induces a map on the level of braid groups. However, note that
this map is neither injective nor surjective. Some vestiges 
of the planar braid groups do however survive:

\begin{thm}
\label{thm_Torsion}
For any tree $\Up$ and any $N>0$, the fundamental group 
$\pi_1(\Cs^N(\Up))$ is torsion-free.
\end{thm}
\pf
A theorem of P. Smith (see \cite[p. 17]{Han91}) 
implies that for a compact 
finite-dimensional Eilenberg-MacLane space $X$ of type 
$(\pi_1(X),1)$, the fundamental group is torsion-free. 
Hence, we may use Corollary~\ref{cor_EM}, noting that one
can deformation retract $\Cs^N(\Up)$ to a compact complex
by enlarging the diagonal slightly.

We give another proof that does not depend upon the 
information concerning higher homotopy groups. The space 
$C^N(\Up)$ is built from products of configuration spaces of 
subgraphs via gluing together along the sets $\Sigma$. These
gluing regions are highly disconnected; however, one may 
perform the requisite gluings in a sequential order. By 
Van Kampen's theorem, the effect on $\pi_1$ of gluing two disjoint spaces
together along a connected set $S$ gives an amalgamated free 
product of the components over $\pi_1(S)$. 
Likewise, gluing a connected set along 
two disjoint copies of a connected subset $S$ yields an 
HNN extension of the fundamental group over $\pi_1(S)$. 

The theorem is trivially true when $N=1$, as well as when there
are no essential vertices. Hence, we may use 
Equation~(\ref{eq_Decomp}) for an induction argument as 
follows. By Lemma~\ref{lem_Injectivity} and the induction 
hypothesis on $N$ and $V$ that fundamental groups of the
pieces are torsion free, we have either an HNN extension or
an amalgamated free product of torsion free groups over an 
injective subgroup. By standard results in geometric group 
theory (see \cite[p.6-8]{Ser80},\cite{SW79}) 
the resulting fundamental group is torsion-free.
\qed

Results about the homotopy type of the configuration spaces are
of interest in understanding the topology of the spaces, but a 
different set of issues must be dealt with in order to have a
practical solution to the control-theoretic problems mentioned in 
the Introduction. For example, the configuration space of $N$ 
points on $\Up$ is an $N$-dimensional complex. Any simplification 
of the space which reduces this dimension will more easily permit 
the construction of explicit control vector fields on the space. 
We offer as a partial solution  to this dilemma a bound on the 
dimension of these spaces up to deformation retraction.

\begin{thm}
\label{thm_DefRetTree}
For any tree $\Up$ and any $N>0$, the configuration space 
$\Cs^N(\Up)$ deformation retracts onto a $V$-dimensional 
subcomplex, where $V$ is the number of essential vertices of $\Up$. 
\end{thm}
\pf 
Choose $p$ and $q$ vertices such that $p\in\del\Up$ and $q$ is 
separated from $p$ by an edge $e$.
We induct on the number of points $N$, the number
of edges $K$ incident to $q$, and the number of essential vertices
$V$ in $\Up$. In order to later apply this argument in the case where
$\Up$ is a general graph, 
the precise induction hypothesis will be that the configuration 
space deformation retracts as a pair $(\Cs^n(\Ga),\Psi)$, where
$\Psi$ denotes the subcomplex 
\be
\Psi := \{\bfx\in\Cs^N(\Up) : x_i\in\del\Up {\mbox{ for some }} i\},
\ee
and $\del\Up$ denotes the boundary of $\Up$. In other words, 
the restriction of the deformation retraction to $\Psi$ induces 
a deformation retraction of $\Psi$. 
If $V=0$ and $N>1$, or if $V>0$ and $N=1$, 
the configuration space pair deformation retracts 
to at most a $V$-dimensional subcomplex and thus satisfies the conclusion. 
Assume, then, that the result holds for all graphs with less than 
$V$ essential vertices. If $K=2$, then the vertex is not essential
and the conclusion is true by induction on $V$. 

Again denote by $\Si_n$ the $n^{th}$ ``end'' of the configuration
space, homeomorphic to $\Cs^{N-1}(\Up)$.
Write $C^N(\Up)$ as the union of pieces 
(as in Equation~(\ref{eq_Decomp}))
\be 
C^N(\Up) \cong C^N(\Up')\bigcup_{\Si'_n}^{n=1..N}\left(
	(\Si_n\times(0,1])\cup(\Si'_n\times\{0\}\right)) .
\ee
Each $\Si_n\times(0,1]$ product end is attached along $\Si'_n
\times\{0\}$; hence, this product deformation retracts 
to $(\Si'_n\times[0,1))\cup(\Si_n\times\{1\})$ rel 
$(\Si'_n\times\{0\})\cup(\Si_n\times\{1\})$ 
as follows. For $x$ not in a neighborhood of $\Si'_n$, shrink 
the segment $\{x\}\times(0,1]$ to $\{x\}\times\{1\}\subset\Si_n$, 
using a bump function with support on a neighborhood of $\Si'_n$. 
Then use another bump function to collapse this 
to $\Si'_n\times\{t\}$ for all $t\in(\eps,1-\eps)$ for $\eps$ small. 
Rounding out the corners in the standard way finishes this stage of 
the deformation. By induction on $V$, the subsets $\Si'_n$ 
deformation retract to subcomplexes of dimension at most $V-1$.
Hence, we may deform each $\Si'_n\times\{t\}$ to a subcomplex of 
dimension at most $(V-1)+1=V$. 

By induction on $K$, one may then deformation retract the base 
$\Cs^N(\Up')$ to a subcomplex of dimension at most $V$ 
without pushing $\Si'$ off of itself. An application 
of the Homotopy Extension Property yields a deformation retraction 
of the entire complex. As a final step, we may 
by induction on $N$ deformation retract the ends 
$\Si_n\times\{1\}$ within themselves and use the Homotopy 
Extension Property again to extend this to a global deformation.
\qed

This deformation retraction is simple enough that one may be 
able to specify a control vector field on the simplified configuration 
space, and then computationally invert the deformation 
retraction to lift this to a control field on the original space. 

% &&&&&&&&&&&&&&&&&&&&&&&&&&&&&&&&&&
\section{Configuration spaces of graphs}
% &&&&&&&&&&&&&&&&&&&&&&&&&&&&&&&&&&

In order to extend the results of the previous section to configuration
spaces of general graphs, some additional techniques and insights
are requisite. Not all results carry over perfectly. For example, 
the configuration space of $N$ points on a circle always 
deformation retracts to a circle, even though this graph has no 
essential vertices. Fortunately, this is the exception and not the 
rule.

In what follows we will denote by $\Ga$ a general (\ie, not 
necessarily simply connected) graph, reserving $\Up$ for trees. 
Choose $P:=\{p_i\}_1^M$ a collection of points in the edge set of $\Ga$ 
such that $\Ga-P$ is a connected open tree. Denote by $\Up$ the 
graph obtained from $\Ga-P$ by adding distinct endpoints; hence, 
every point $p_i\in P$ is split into two points $p_i^+$ and $p_i^-$ in $\Up$. 
The configuration space $\Cs^N(\Ga)$ decomposes as $\Cs^N(\Up)$ 
with certain ends identified pairwise.
More specifically, the graph $\Up$ has $2M$ ``ends'' corresponding to 
the points $p_i^\pm$. Likewise, for each such end of $\Up$, the 
configuration space $\Cs^N(\Up)$ has $N$ ``ends'' which come from 
Equation~(\ref{eq_Decomp}). Hence, the configuration space 
$\Cs^N(\Up)$ has $2MN$ product ends.

\begin{thm}\label{cor_EMGraph}
The configuration space $\Cs^N(\Ga)$ is a $K(\pi_1,1)$.
\end{thm}
\pf
Decompose the configuration space by splitting along the 
aforementioned set $P$. Then, given a representative $f:S^k\ra
\Cs^N(\Ga)$ of $\pi_k$, we know that $f$ must intersect the 
clipped ends nontrivially via Corollary~\ref{cor_EM}. However, 
the proof of Lemma~\ref{lem_Injectivity} is presented in the 
context of contracting a loop in a product end of a tree; 
hence, the inclusion of the ends of $\Cs^N(\Ga)$ into 
$\Cs^N(\Up)$ is $\pi_1$-injective, and the proof of 
Theorem~\ref{thm_Contract} applies.
\qed

\begin{cor}\label{cor_TorFreeGraph}
The pure braid group of a graph is torsion-free.
\end{cor}

\begin{thm}\label{thm_DefRetGraph}
For any graph $\Ga$ not homeomorphic to a circle, the configuration 
space $\Cs^N(\Ga)$ deformation retracts to a $V$-dimensional 
subcomplex, where $V$ is the number of essential vertices in $\Ga$. 
\end{thm}
\pf
In the deformation retraction of Theorem~\ref{thm_DefRetTree}, 
the deformation can always be accomplished rel the end $\Psi$ (except 
when the tree is the trivial line segment with one point on it).
Hence, we may deformation retract that portion of $\Cs^N(\Ga)$ 
which corresponds to $\Cs^N(\Ga-P)$ down to a $V$-dimensional
subcomplex. Then, the remaining portions of the configuration 
space may be likewise deformation retracted by induction on 
the number of points as in the previous theorem. 
\qed

We close our treatment without a specific classification of 
the configuration spaces of graphs. By Corollary~\ref{cor_EMGraph}, 
the isomorphism class of the pure braid groups on graphs 
determines the homotopy type of the configuration space. Given 
the examples from the next section, we conjecture that the 
configuration spaces are all homotopy equivalent to collections of tori of 
various dimensions, glued together along incompressible 
subtori: a reasonable refinement of this would be the following.
\begin{conj}
\label{conj_Artin}
The group $\pi_1(\Cs^N(\Up))$ is an {\em Artin right angle group}: 
the group has a presentation in which all of the relations are 
commutators of the generators. 
\end{conj}
%

% &&&&&&&&&&&&&&&&&&&&&&&&&&&&&&&&&&
\section{Examples}
% &&&&&&&&&&&&&&&&&&&&&&&&&&&&&&&&&&

{\sc example 1:} Let $\Cs^N_K$ denote the configuration space
of $N$ points on a $K$-pronged radial tree (\ie, having $K$
edges attached to a single central vertex). According to 
Theorem~\ref{thm_DefRetTree}, $\Cs^N_K$ deformation retracts to 
a one-dimensional graph. Since the homotopy type of a graph 
is determined by its Euler characteristic, we can derive the 
following:
\begin{prop}
\label{prop_CNK}
The braid group $\pi_1(\Cs^N_K)$ is isomorphic to a free group
on $Q$ generators, where $Q$ equals
\be
	Q := 1 + \left(NK-2N-K+1\right)\frac{(N+K-2)!}{(K-1)!}
\ee
\end{prop}
\pf
Using Equation~(\ref{eq_Decomp}), one derives a recursion 
relation for the Euler characteristic:
\be
\label{eq_Euler}
\chi(\Cs^N_K) = \chi(\Cs^N_{K-1})+ N\left[\chi(\Cs^{N-1}_K)
		- E\right] ,
\ee
where $E$ denotes the number of connected components of 
$\Si-\iota(\Si')$.
Each such component contributes one edge in the 
deformation retracted space and hence contributes a $-1$
to the value of $\chi(\Cs^N_K)$. A simple combinatorial 
argument shows that 
\be
	E = \prod_{i=1}^{K-1}(N+i-2) = \frac{(N+K-3)!}{(K-2)!} .
\ee
The seed for the recursion relation is the fact that 
$\Cs^1_K$ is homeomorphic to the underlying tree which is 
homotopically trivial. Solving (\ref{eq_Euler}) yields 
\be
	\chi(\Cs^N_K) = -\left(NK-2N-K+1\right)\frac{(N+K-2)!}{(K-1)!},
\ee
which, in turn, implies that the configuration space is 
homotopic to a wedge of $1-\chi$ circles.
\qed

{\sc example 2:} The spaces $\Cs^N_K$ of Example 1 are, by 
Proposition~\ref{prop_CNK}, homotopic to a wedge of circles.
However, the deformation retraction of Theorem~\ref{thm_DefRetTree} 
does not compress the space to this extreme, but rather leaves some 
structure. In the simple case of $\Cs^2_3$, the deformation 
retraction yields a 1-d graph which resembles a `benzene ring': 
homeomorphic to $S^1$ with six radial edges attached. This 
reduction of the configuration space has the advantage that 
each vertex corresponds to a necessary passage through the 
central vertex. 

% =========================================
\begin{figure}[htb]
\begin{center}
	\epsfxsize=5.0in\leavevmode\epsfbox{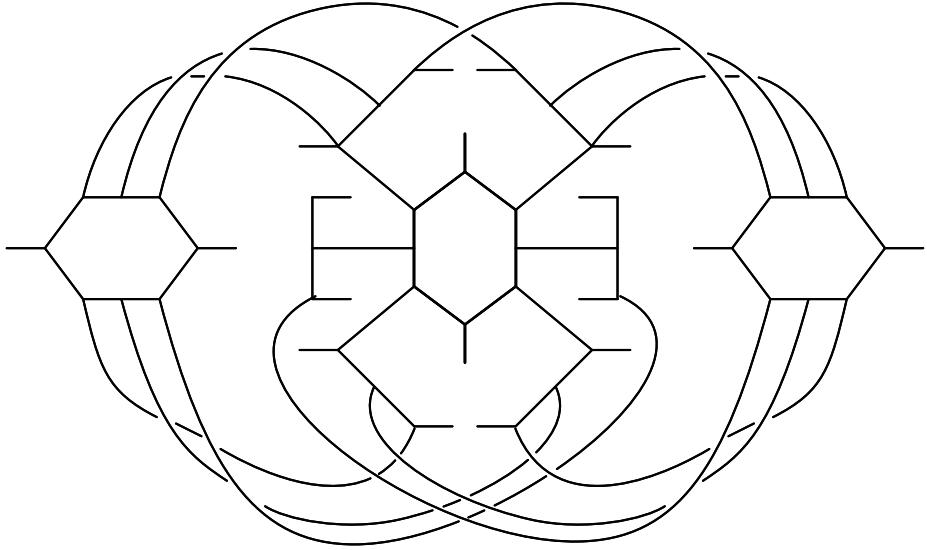}
\end{center}
\caption{The space $\Cs^3_3$ deformation retracts to a graph.} 
\label{fig_Hanoi}
\end{figure}
% =========================================

In Figure~\ref{fig_Hanoi}, we illustrate this semi-minimal reduction 
of the space $\Cs^3_3$. According to Proposition~\ref{prop_CNK}, 
$\Cs^3_3$ deformation retracts to the wedge of 13 circles. However, 
under the piecewise deformation given by the proof of 
Theorem~\ref{thm_DefRetTree}, one obtains a graph that is 
non-planar and clearly has vestigial copies of the `benzene' graph 
of $\Cs^2_3$. It is interesting to note that this is the configuration
space associated to a type of ``Towers of Hanoi'' problem:
the condition that one may only move one ring at a time
corresponds to the condition that the central vertex may be 
occupied by at most one point. The diameter of the graph in 
Figure~\ref{fig_Hanoi} corresponds to the minimal number of
steps required to reverse the order of three points initially on the
same edge. 

{\sc Example 3:} Consider the graph $\Ga_H$ which is 
homeomorphic to the letter ``H'': two essential vertices. 
According to Theorem~\ref{thm_DefRetTree}, the configuration 
space $\Cs^N(\Ga_H)$ deformation retract to a 2-complex. 
However, a more careful analysis of individual cases 
yields more insightful results. In what follows, we have 
executed the proof of Theorem~\ref{thm_DefRetTree} step-by-step, 
applying the results of Example 1 at each stage. We 
denote by $F_p$ the free group on $p$ generators.

{\sc example 3a:} $\pi_1(\Cs^2(\Ga_H))\cong F_{3}$.

{\sc example 3b:} $\pi_1(\Cs^3(\Ga_H))\cong F_{25}$.

{\sc example 3c:} $\pi_1(\Cs^4(\Ga_H))\cong 
%	F_{189}*\left[*_{1}^6\left((\zed\times\zed)*\zed\right)\right]$.
	F_{195}*\left[*_{1}^6\left(\zed\times\zed\right)\right]$.

	There are exactly six 2-tori in a homotopically minimal
	representative of this space. Each torus corresponds to 
	a configuration where pairs of points trace out loops
	in a neighborhood of the individual vertices.

It is entirely possible that Example 3c illustrates the canonical way in 
which non-free components of braid groups on graphs can arise:
\cf Conjecture~\ref{conj_Artin}. 

\begin{figure}[htb]
\begin{center}
	\epsfxsize=2.75in\leavevmode\epsfbox{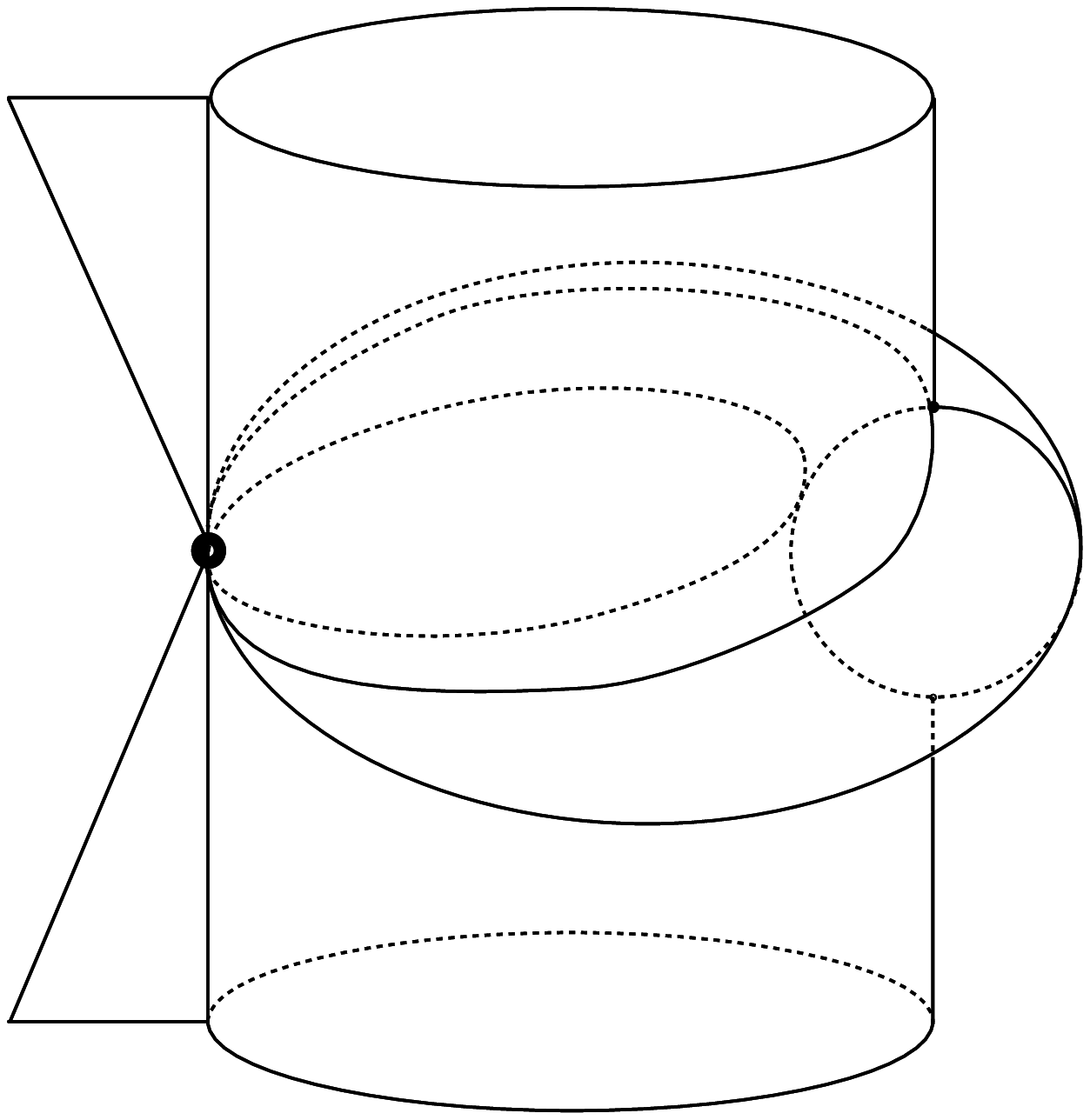}
\end{center}
\caption{The space $\Cs^2(\Ga_Q)$ embeds in $\real^3$.} 
\label{fig_Q}
\end{figure}

{\sc example 4:} The simplest non-tree graph, 
$\Ga_Q$, is homeomorphic to a circle with one edge attached.
As this graph is the identification of two edges of the 
Y-graph $\Ga_Y$, one may obtain the configuration space 
(up to homeomorphism) via the proper identifications. We
display the result in Figure~\ref{fig_Q}. Note that this 
space deformation retracts onto a one-dimensional graph
homotopic to the wedge of three circles (the point on the 
left of the diagram is a puncture). 

% =========================================
% =========================================

% ACKNOWLEDGMENT

\small

{\sc acknowledgment}: This paper was inspired via the vision 
and enthusiasm of Dan Koditschek. Conversations with Mladen 
Bestvina, John Etnyre, Emily Hamilton, Allen Hatcher, and Alec Norton were 
of great aid. This work was supported in part by the National 
Science Foundation [DMS-9508846].

% BIBLIOGRAPHY
\bibliographystyle{alpha}

\end{document}